\documentclass[12pt]{article}
\usepackage[centertags]{amsmath}
\usepackage{amsfonts}
\usepackage{amssymb}
\usepackage{amsthm}
\usepackage{amsmath,mathrsfs}
\usepackage{amsmath,amscd}
\usepackage{mathtools}
\usepackage[matrix,arrow,curve,cmtip]{xy}
\title{\textbf{The Physical Mathematics of Segal Topoi and Strings}}
\author{Renaud Gauthier \footnote{2020 Math. Subj. Class: 18N60, 14D23, 18F20, 18A22, 81T30  Keywords: Segal topoi, derived stacks, quantum states, dynamics, strings, M-theory.} \\ \\}
\theoremstyle{definition}

\DeclareMathOperator*{\colim}{\text{colim}}

\newcommand{\beq}{\begin{equation}}
\newcommand{\eeq}{\end{equation}}

\newcommand{\rarr}{\rightarrow}

\newcommand{\Rarr}{\Rightarrow}

\newcommand{\Ob}{\text{Ob\,}}

\newcommand{\xrarr}{\xrightarrow}

\newcommand{\cA}{\mathcal{A}}

\newcommand{\cC}{\mathcal{C}}

\newcommand{\cD}{\mathcal{D}}

\newcommand{\cO}{\mathcal{O}}

\newcommand{\cX}{\mathcal{X}}

\newcommand{\bL}{\mathbb{L}}

\newcommand{\bR}{\mathbb{R}}

\newcommand{\bra}{\langle}

\newcommand{\Fun}{\text{Fun}}

\newcommand{\Hom}{\text{Hom}}
\newcommand{\Ho}{\text{Ho}\,}

\newcommand{\ket}{\rangle}

\newcommand{\Mor}{\text{Mor}\,}
\newcommand{\Map}{\text{Map}}

\newcommand{\Nat}{\text{Nat}}
\newcommand{\op}{\text{op}}

\newcommand{\Set}{\text{Set}}
\newcommand{\Spec}{\text{Spec\,}}
\newcommand{\Top}{\text{Top}}
\newcommand{\uHom}{\underline{\Hom}}

\newcommand{\ACommC}{A\text{-Comm}(\cC)}

\newcommand{\AMod}{A\text{-Mod}}

\newcommand{\AffC}{\cA \text{ff}_{\cC}}
\newcommand{\AffChat}{\AffC^{\,\wedge}}
\newcommand{\AffCtildetau}{\AffC^{\,\sim,\tau}}

\newcommand{\Comm}{\text{Comm}}

\newcommand{\CommC}{\text{Comm}(\cC)}

\newcommand{\Der}{\mathbb{D}\text{er}}

\newcommand{\dkAff}{\text{d}k\text{-Aff}}

\newcommand{\dStk}{\text{dSt}(k)}

\newcommand{\DkAff}{\text{D}^- k \text{-Aff}}

\newcommand{\eps}{\epsilon}

\newcommand{\ffqc}{\text{ffqc}}

\newcommand{\LFx}{\bL_{F,x}}

\newcommand{\oP}{\oplus}

\newcommand{\oT}{\otimes}
\newcommand{\oTA}{\oT_A}
\newcommand{\oTAL}{\oTA^{\mathbb{L}}}

\newcommand{\RHom}{\mathbb{R} \uHom}

\newcommand{\RHomLg}{\RHom^*_{\SeT}}

\newcommand{\RuHom}{\bR \uHom}
\newcommand{\RuSpec}{\bR \uSpec}

\newcommand{\sPr}{\text{sPr}}

\newcommand{\SetD}{\Set_{\Delta}}

\newcommand{\St}{\text{St}}

\newcommand{\Sp}{\text{Sp}}

\newcommand{\skCAlg}{\text{s}k\text{-CAlg}}

\newcommand{\SeT}{\text{SeT}}

\newcommand{\SeCat}{\text{SeCat}}

\newcommand{\TFx}{\mathbb{T}_{F,x}}

\newcommand{\uh}{\underline{h}}

\newcommand{\uSpec}{\underline{\Spec}}

\newcommand{\dAffC}{\text{d}\AffC}
\newcommand{\dAffChat}{\text{d}\AffChat}
\newcommand{\dAffCtildetau}{\text{d}\AffCtildetau}
\newcommand{\RHomXX}{\RHom(\cX,\cX)}
\newcommand{\Cusp}{\text{Cusp}}
\newcommand{\beeT}{\mathbb{T}}
\newcommand{\TPsix}{\beeT_{\Psi,x}}
\newcommand{\LPsix}{\bL_{\Psi,x}}
\newcommand{\maxMapgen}{\text{max}\Map^{gen}}

\newcommand{\uF}{\underline{F}}

\newcommand{\uu}{\underline{u}}
\newcommand{\utau}{\underline{\tau}}
\newcommand{\utaust}{|\utau \ket}
\newcommand{\upsi}{\underline{\psi}}
\newcommand{\uG}{\underline{G}}
\newcommand{\uGst}{|\uG \ket}
\newcommand{\uFst}{|\uF \ket}
\newcommand{\dX}{\delta \cX}
\newcommand{\dtwoX}{\delta^2 \cX}
\newcommand{\dnX}{\delta^n \cX}
\newcommand{\dmX}{\delta^m \cX}

\newcommand{\psinm}{\psi^{(n,m)}}
\newcommand{\End}{\text{End}}
\newcommand{\dXgen}{\dX^{gen}}
\newcommand{\cCgen}{\cC^{gen}}
\newcommand{\XoTDS}{\cX \oT^{\Delta} \Sigma}

\begin{document}
\maketitle
\begin{abstract}
	We introduce a notion of dynamics in the setting of Segal topos, by considering the Segal category of stacks $\cX = \dAffCtildetau$ on a Segal category $\dAffC = L \CommC^{\op}$ as our system, and by regarding objects of $\RHomXX$ as its states. We develop the notion of quantum state in this setting and construct local and global flows of such states. In this formalism, strings are given by equivalences between elements of commutative monoids in $\cC$, a base symmetric monoidal model category. The connection with standard string theory is made, and with M-theory in particular.  
\end{abstract}

\newpage

\section{Introduction}
The present paper has two aims, the first one being to develop a notion of dynamics in towers of Segal topoi, any Segal topos in the tower being seen as embedded in the one that follows, and the second aim being an application of this formalism to Physics, and to strings in particular, all of this being developed in the spirit of Physical Mathematics in the sense Greg Moore gives it (\cite{GM1}, \cite{GM2}, \cite{BFMNRS}).\\

It was in \cite{RG} that we focused on the Segal category of geometric morphisms $\RHomLg(\cX,\cX)$, with $\cX = \dStk = L(\DkAff^{\, \sim, \ffqc})$ (\cite{T}), for the sake of studying Higher Galois. Recall that $\dStk = \dkAff^{\, \sim, ffqc}$ is the localization of $\widehat{\dkAff} = \RHom(\dkAff^{\op}, \Top)$ with respect to ffqc hypercovers, where $\Top = L \SetD$ (\cite{T}, \cite{TV1}). Such a Segal category is a Segal groupoid however, and with a view towards applying the present formalism to Physics, where the notion of entropy is fundamental, it is therefore natural to just focus on $\RHomXX$ instead. We will work in full generality, and take a base Segal category built from a generic (simplicial) symmetric monoidal model category $\cC$ satisfying the conditions of a homotopical algebraic context as enunciated in \cite{TV4}. Under these conditions, $\dAffC = L \CommC^{\op}$ is our base Segal category on which we put a Segal topology $\tau$, giving rise to a Segal category of (derived) stacks $\dAffCtildetau$, a Segal topos. For the sake of applications to Physical Mathematics, this is really the Segal topos we have in mind, but for a formal development of the underlying theory, we will work with a generic Segal topos $\cX$.\\

Identifying $\cX$ with $id_{\cX} \in \RuHom(\cX,\cX)$, we regard objects of $\RHomXX$ other than the identity as providing incremental changes in our system $\cX$. Indeed, an instantaneous snapshot of our system is given by providing all objects of $\cX$ and the morphisms between them. For instance, a morphism $f:F \rarr G$ in $\cX$ relates $F$ and $G$ together. Now given $\zeta \neq id$ in $\RHomXX$, $\zeta F \rarr \zeta G$ provides a modification to such a relation, and this for all such morphisms $f$, and for all objects $F$ and $G$ of $\cX$. It is because such functors $\zeta$ provide an overall change in our system $\cX$ that we regard them as states of our system. For that reason we will refer to $\RHomXX$ as the space of states of our system $\cX$. \\

\newpage
$F \in \cX$ being fixed, if one initially obtains a shift $F \mapsto \zeta F$ upon applying $\zeta$, using $\zeta$ again one obtains a sequence of incremental changes in the dynamics of $F$ as in $F \mapsto \zeta F \mapsto \zeta^2 F \mapsto \cdots$. The colimit $\zeta^{\infty} F = \colim \zeta^n F$ of this sequence in $\cX \in \SeT$ (the Segal category of Segal topos) provides a discrete flow of $F$ in the direction of $\zeta$. A (smooth) flow is obtained if $\zeta$ is an infinitesimal state. In order to characterize such states, we will develop a notion of quantum state, being defined as any functor $\zeta: \cX \rarr \cX$ that is minimal in the way it modifies objects of $\cX$, something that we will define formally in Section \ref{qs}. Now, there is no reason for keeping $\zeta$ throughout. Thus we consider how such functors $\zeta$ evolve, and if we do this functorially we are led to considering functors $\Psi: \RHomXX \rarr \RHomXX$. Iteratively applying $\Psi$ to $\zeta$, we obtain a sequence:
\beq
\zeta \rarr \Psi \zeta \rarr \Psi^2 \zeta \rarr \cdots \nonumber 
\eeq
with a colimit $\Psi^{\infty} \zeta = \colim \Psi^n \zeta$, which once applied to $F \in \cX$ provides a discrete flow
\beq
F \rarr \zeta F \rarr \Psi \zeta F \rarr \Psi^2 \zeta F \rarr \cdots \nonumber
\eeq
for $F$, which becomes a smooth, physical flow if both $\zeta$ and $\Psi$ are quantum states. Iterating this construct, there is no reason that $\Psi$ always be used to modify $\zeta$, so we should consider how $\Psi$ itself evolves. Adopting the notation $\RHomXX = \dX$, we have $\Psi \in \dtwoX$, and we are therefore looking for an overarching object $\Gamma \in \delta^3 \cX$ to repeat this construction. Formalizing this, if $\psi^{(-,m)} \in \delta^m \cX$ and $\psi^{(n+1,m)} = \psi^{(n,m+1)} \psinm$, we get a discrete flow at the level of $\cX$:
\beq
F \xrarr{\psi^{(0,1)}} \psi^{(0,1)} F \xrarr{\psi^{(1,1)}} \psi^{(1,1)} \psi^{(0,1)}F \xrarr{\psi^{(2,1)}} \psi^{(2,1)} \psi^{(1,1)} \psi^{(0,1)} F \rarr \cdots \nonumber
\eeq
where the objects $\psi^{(n,1)}$ of $\dX$ fit into a sequence:
\beq
\psi^{(0,1)} \xrarr{\psi^{(0,2)}} \psi^{(0,2)} \psi^{(0,1)} = \psi^{(1,1)} \xrarr{\psi^{(1,2)}} \psi^{(1,2)} \psi^{(1,1)} = \psi^{(2,1)} \rarr \cdots \nonumber
\eeq
and the objects $\psi^{(n,2)}$ of $\dtwoX$ themselves fit into a sequence:
\beq
\psi^{(0,2)} \xrarr{\psi^{(0,3)}} \psi^{(0,3)} \psi^{(0,2)} = \psi^{(1,2)} \xrarr{\psi^{(1,3)}} \psi^{(1,3)} \psi^{(1,2)} = \psi^{(2,2)} \rarr \cdots \nonumber
\eeq
ad infinitum. We argue $\colim_n \dnX$ contains all the possible higher dynamics of objects of $\cX$, which exists by virtue of \cite{HS}, or even \cite{RV}. Further, the flow of some $F \in \cX$ can be made smooth by taking all the higher states $\psinm$ involved to be quantum in nature, for all $m \geq 1$ and all $n \geq 0$.\\

For simplicity, and for illustrative purposes, we will focus on the case where we have one single state $\zeta = \psi^{(0,1)} \in \dX$, one $\Psi = \psi^{(0,2)} \in \dtwoX$, and we consider a (local) flow of $\zeta$ induced by $\Psi$. We are led to introducing the notion of cusp at $(\Psi, \zeta)$, which we will define by way of a tangent object, using as impetus the derivation of the cotangent complex the way it was developed in \cite{TV4}. To be precise, $\Cusp[\Psi, \zeta] = \{ \Phi: \dX \rarr \dX \, | \, \beeT_{\Psi, \zeta} \cong \beeT_{\Phi,\zeta} \}$, where $\beeT_{\Psi, \zeta}$ is the tangent object of $\Psi$ at $\zeta$, defined as the dual of the cotangent object $\bL_{\Psi, \zeta}$, a crude first approximation thereof being provided by $ \Map (\Psi \zeta, \Psi \rho \circ \Psi \zeta)$ for some $\rho \in \dX$. We did not yet specify the exact nature of such objects, the reason being that we will give two presentations of this formalism, one being Set-based, the other one being simplicial in nature, and accordingly, depending on the context, one will talk about a (co)tangent set, or a (co)tangent simplex.\\

Finally we apply this formalism to strings by focusing on $\cX = \dAffCtildetau$, which we introduced in the second paragraph above, and by defining a string between $a \in A$, $a' \in A'$ in $\CommC$ to be an equivalence relation defined by $a \stackrel{\gamma}{\sim} a'$ if $A \cong A'$ in $L \CommC$. It then becomes apparent that in algebraic geometry we obtain an algebraic version of string theory. From this perspective, stacks $F: \CommC \rarr \SetD$ map algebras $A \in \CommC$, viewed as natural laws, to simplicial sets $FA$, and $a \stackrel{\gamma}{\sim} a'$ in $A$ maps to $Fa \stackrel{F \gamma}{\sim} Fa'$, a string in $FA$. We regard $\gamma$'s as formal, algebraic strings with realization $F \gamma$ in an ambient space $FA$. All such images $FA$ are in a sense representations of a law $A$, and are related to one another by functoriality of $F$. Since we expect physical phenomena to be derived from the dynamics of strings, and to exhibit natural laws in a coherent fashion, we are therefore led to viewing such functors $F$ as modeling physical phenomena, some of which in particular will be physical fields.\\

In section \ref{base Segal topos} we fix notations about our base symmetric monoidal model category $\cC$ and construct the Segal category of stacks $\dAffCtildetau$. In section \ref{Dynamics}, we start with a Segal topos $\cX$, which we consider as our system, and consider states on it, morphisms of states, and accompanying flows. The iterative picture is presented and discussed. We also provide an analysis of what would make a state quantum in nature for the sake of having a working definition of flows. In section \ref{Cusps} we focus on local flows. In order to do this we have to introduce the concepts of generalized categories, and of tangent and cotangent objects. Finally we apply this formalism to the study of strings in the last section, which we regard as equivalences between objects of commutative monoids of $\cC$, which was fixed at the onset. In particular we explain why the present setting is a natural place to have strings; not only does it explain that strings are not seen in nature, but it also suggests that a functorial treatment may be a good place to develop M-theory.\\

Notations: we denote by $\SetD$ the category of simplicial sets, by $\SeCat$ the category of Segal categories, and by $\SeT$ the Segal category of Segal topos. Often we will abbreviate $\RHomXX$ by just $\dX$.

\section{Building up the base Segal topos} \label{base Segal topos}
As mentioned in the introduction, in lieu of working with $\dkAff = L \skCAlg^{\op}$ the way we would in derived algebraic geometry, we abstract this formalism and start with a base (simplicial) symmetric monoidal category $\cC$ endowed with a model category structure compatible with its monoidal structure. This makes $\cC$ into a symmetric monoidal model category. We ask that it be combinatorial for localization purposes, that it be pointed, that $\Ho(\cC)$ be additive, and that $\Comm(\cC)$ as well as modules over its objects have a model category structure. Those conditions are covered in \cite{TV4} and collectively they define a homotopical algebra (HA) context. To be precise, asking that $\cC$ be pointed and combinatorial, and that $\Ho(\cC)$ be additive, this corresponds to Assumption 1.1.0.1 of \cite{TV4}. For later purposes, if $A \in \Comm(\cC)$, we want $\AMod$ to be combinatorial, proper, with its induced tensor product $- \oT_A -$ making it into a symmetric monoidal model category. This is Assumption 1.1.0.2 of the same reference. Assumption 1.1.0.3 ensures that $- \oT_A M$ preserves equivalences. Assumption 1.1.0.4 ensures that $\ACommC$ as well as $(\ACommC)_{nu}$ ($nu$ stands for non-unital) are combinatorial, proper, and if $B$ is cofibrant in $\ACommC$, then $B \oT_A -$ preserves equivalences. Now less obvious is the following assumption which we will adopt for the sake of infinitesimal lifting properties: if we fix a subcategory $\cC_0$ of $\cC$, Assumption 1.1.0.6 asks that $1 \in \cC_0$, that $\cC_0$ be stable by equivalences and homotopy colimits, and that $\Ho(\cC_0) \subset \Ho(\cC)$ be stable by the derived tensor product. \\

All this being defined, and using the notations and definitions of \cite{TV4},  our base model category $\Comm(\cC)^{op}$ is denoted $\AffC $, whose objects are referred to as affine objects (derived affine objects if $\cC$ is simplicial). 
If $A \in \Comm(\cC)$, the corresponding object in $\AffC$ is denoted $\Spec A$. The category of (derived) prestacks on $\AffC$ is denoted $\AffChat$, it is the left Bousfield localization of the category of simplicial presheaves $\sPr(\AffC)$ with respect to equivalences, that is $\AffChat = L_{Bous} \sPr(\AffC)$. Its objects are functors $F: \AffC^{op} \rarr \SetD$ that preserve equivalences. One denotes by $\RHom$ the derived simplicial Hom of the simplicial model category $\AffChat$. If we fix a model topology $\tau$ on $\AffC$, one obtains a model category of stacks $\AffCtildetau$ on $(\AffC, \tau)$, whose homotopy category $\St(\cC, \tau) = \Ho(\AffCtildetau)$ is the category of (derived) stacks on $\AffC$ for the topology $\tau$. Its objects are functors $F: \AffC^{\op} \rarr \SetD$ that preserve equivalences, and satisfy the descent property for homotopy $\tau$-hypercovers. One assumption we need at this point is Assumption 1.3.2.2 of \cite{TV4} that asks that $\tau$ on $\Ho(\AffC)$ be quasi-compact, plus two other technical requirements that the reader can find in \cite{TV4}, both of which being needed for the functorial behavior of representable stacks under equivalences, as well as for the convenient characterization of $\AffCtildetau$ as the left Bousfield localization of $\AffChat$. This construction of $\AffCtildetau$, starting from a model category $\AffC$, can equivalently be done starting from a simplicial site $(T,\tau)$ instead, with $T = L \AffC$, giving rise to a model category of stacks $\sPr_{\tau}(T)$.\\

All that was done so far does not yet involve Segal categories. The transition is being made by regarding the simplicial category $\dAffC = L(\CommC^{\op})$ as an object of $ \Ho(\SeCat)$, whose objects are generically referred to as (derived) affine schemes, following standard usage (\cite{T}). From \cite{TV1}, one defines $\dAffChat = \RHom(\dAffC^{\op}, \Top)$ - where $\Top = L\SetD$ - referred to as the Segal category of (derived) prestacks on $\dAffC$. Let $\tau$ be a Segal topology on $\dAffC$, that is a Grothendieck topology on $\Ho(\dAffC)$.  Then the Segal category of (derived) stacks $\dAffCtildetau$ on $(\dAffC, \tau)$ is a full subcategory of $\dAffChat$ (a Segal topos as well), that one can construct as a left exact localization of $\dAffChat$, or as a left Bousfield localization thereof with respect to $\tau$-hypercovers. \\

\newpage

The connection with our initial construction starting from a simplicial site $(T,\tau)$ is given by:
\beq
L(\sPr_{\tau}(\dAffC)) \cong \dAffCtildetau \nonumber
\eeq
where $L$ is a simplicial localization of model categories, we used the fact that $\dAffC = L(\Comm(\cC)^{\op})$ is a simplicial category, and the fact that a topology on the simplicial category $\dAffC$ also corresponds to a Segal topology $\tau$ on $\dAffC$ viewed as a Segal category, since simplicial sites and Segal sites are defined similarly (\cite{TV1}, \cite{TV4}).\\

\section{Dynamics} \label{Dynamics}
We regard $\cX = \dAffCtildetau$ as modeling natural phenomena, for an appropriately chosen fundamental (simplicial) symmetric monoidal model category $\cC$ and a topology $\tau$ on $\dAffC$, both providing a homotopical algebraic geometry context (see \cite{TV4} for details). In this picture, natural laws are modeled by objects of $\CommC$, and a given phenomenon in this picture is represented by a stack $F \in \cX$ valued in $\SetD$, since $\cX$ can be modeled by functors $F: \Comm(\cC) \rarr \SetD$. The problem is that this is very much inert; for $A \in \CommC$, $F(A) \in \SetD$ presents how $F$ responds to a law $A$. For another $B \in \CommC$, assuming we have a morphism $f: A \rarr B$ between those laws, by functoriality we have a map $Ff: F(A) \rarr F(B)$, thereby displaying some coherence in the representations of $F$, but there is no manifest dynamics. In other terms $F \in \cX$ provides one coherent manifestation of physical laws, and morphisms in $\cX$ relate different manifestations together, without providing any sense in which those stacks are deformed to produce dynamics.\\

In what follows, $\cX \in \SeT$. We identify $\cX$ with $id_{\cX} \in \RHomXX$, the identity therefore providing a snapshot of inter-relations between objects of $\cX$ as argued above, so that other non-identity morphisms $\psi: \cX \rarr \cX$ in $\RHomXX$ can consequently be interpreted as providing an overall shift in the objects $F$ of $\cX$, a first step towards defining a dynamics in $\cX$. Yet any such morphism $\psi$ is but an instantaneous deformation; $\psi: \cX \rarr \cX$ maps any $F \in \cX$ to $\psi F$. For a smoother representation, since this is just incremental, we seek a flow generated by such functors $\psi$. The most basic flow consists in having $\psi \in \RHomXX$ fixed and on repeatedly applying it to obtain a sequence of the form:
\beq
F \xrarr{\psi} \psi F \xrarr{\psi} \psi^2 F \xrarr{\psi} \cdots \nonumber
\eeq
More generally, one can consider various objects $\psi = \psi_0, \psi_1, \cdots $ of $\RHomXX$ being used to construct a flow, as in:
\beq
F \xrarr{\psi_0} \psi_0 F \xrarr{ \psi_1} \psi_1 \psi_0 F \xrarr{\psi_2} \psi_2 \psi_1 \psi_0 F \xrarr \cdots \nonumber
\eeq
Such a flow can be provided by a functor $\Psi: \RHomXX \rarr \RHomXX$ for which $\Psi \psi_i = \psi_{i+1}$. Those are but the first stages of the iterative construction presented in the introduction. In the present section we motivate referring to objects of $\RHomXX$ as states, we construct flows, and we study the iterative picture.\\

\subsection{States and operators}
\subsubsection{General formalism}
For simplicity we write $\dX$ for $\RHomXX$. We informally refer to $\cX$ as a \textbf{system}, and we define a \textbf{state} of a system $\cX$ as being an object of $\dX$, that is a functor of Segal categories $\cX \rarr \cX$. Indeed, $\psi \in \dX$ maps all $F's \in \cX$ to $\psi F \in \cX$, and therefore corresponds in this picture to an overall change in the system $\cX$, whence the definition of state. We can formalize this as follows. Amongst objects $\psi \in \dX$, there are those functors that change the intrinsic nature of objects $F \in \cX$, and those that merely casts them in another state without modifying them otherwise. We are interested in the latter case of transformation. If $\psi F$ is $F$ shifted to a certain state, in the physical sense of the word, this implies $F$ itself was initially in a given state. This means we are not modifying the intrinsic object $F$, but the state it finds itself in. To make this manifest, we introduce $\uF$ the underlying object associated with $F$ prior to being put in a certain state, it being understood that the latter is built from $\uF$ and a certain physical state. We introduce a symmetric monoidal category of (physical) states $\Sigma$, whose objects are denoted by $|-\ket$, which we assume form a complete, orthonormal set, with product $m: \bra p| \oT |q \ket \rarr \bra p | q \ket = \delta_{pq}$. We regard $\uF$ as an object of the same type as $F$, and an operator as well, and if the letter $\cO$ designates operators, we denote by $\cX \wedge \cO$ the category of underlying objects of $\cX$ that are also operators. $\uF$ is an operator insofar as its eigenvalues on $\Sigma$ are objects of $\cX$ as in:
\beq
\uF|n \ket = F_n | n \ket \nonumber
\eeq
where $F_n \in \cX$. In particular we will designate by $|0 \ket$ the state of $\Sigma$ that corresponds to $F$, in such a manner that $\uF|0 \ket = F |0 \ket$. With these notations:
\beq
F = \bra 0 | \uF |0 \ket \nonumber
\eeq
appears as the expectation value of $\uF$. More precisely, we first have operator-states evaluation maps:
\begin{align}
	ev: \cO \oT \Sigma & \rarr \Sigma \nonumber \\
	\Gamma \oT |p \ket & \mapsto \Gamma |p\ket = \sum c_n |n \ket \nonumber
\end{align}
for certain coefficients $c_n$, that one can generalize to:
\begin{align}
	ev: \cX \wedge \cO \oT \Sigma & \rarr \cX \oT \Sigma \label{FOS} \\
	\uF \oT |p \ket & \mapsto \uF| p \ket = F_p |p\ket \nonumber
\end{align}
To be precise, if one defines $\pi: \cX \rarr \Sigma$, $F_p \mapsto |p \ket$, then denoting by $\XoTDS$ the image of the composite $\cX \xrarr{\Delta} \cX \oT \cX \xrarr{id \oT \pi} \cX \oT \Sigma$, the evaluation map in \eqref{FOS} is valued in $\XoTDS$. Focusing now on the product in $\Sigma$, $m: \Sigma \oT \Sigma \rarr \bR$, it has an obvious generalization:
\begin{align}
	m: \Sigma \oT (\XoTDS) & \rarr \cX \label{SFS} \\
	\bra p | \oT \sum F_n |n \ket & \mapsto \sum F_n \bra p | n \ket = F_p \nonumber
\end{align}
Combining \eqref{FOS} and \eqref{SFS}, we get:
\begin{align}
	\Sigma \oT ( \cX \wedge \cO \oT \Sigma) & \xrarr{1 \oT ev} \Sigma \oT (\XoTDS) \xrarr{m} \cX \nonumber \\
	\bra 0 | \oT (\uF \oT |0 \ket) & \mapsto \bra 0 | \oT F|0 \ket \longmapsto \bra 0 | F | 0 \ket = F \nonumber
\end{align}

Alternatively, one may want a state-like representation of $\uF$ for the purpose of doing such evaluations. If $\uF$ is fully defined by its eigenvalues $F_n$ on states $|n \ket$, then it suffices to consider the isomorphism:
\begin{align}
	\cX \wedge \cO & \xrarr{\cong} \XoTDS \nonumber \\
	\uF & \mapsto \uFst = \sum |n \ket \bra n | \uF |n \ket = \sum F_n |n \ket \nonumber
\end{align}
With these notations:
\beq
\bra 0 | \uF \ket = \bra 0 | \sum |n \ket \bra n | \uF |n \ket = \bra 0 | \uF | 0 \ket = F \nonumber
\eeq

\subsubsection{Objects of $\cX$ as functors}
Since $\cX$ is a functor category in practical applications, we focus on this case presently. We want to make sure the functorial calculus at the level of  $\uF \in \cX \wedge \cO$ projects down to the same calculus for objects of $\cX$ in any state. Consider a map between objects of $\cX \wedge \cO$: $\utau: \uF \Rarr \uG$, defined levelwise by $\utau_n: F_n \Rarr G_n$ by virtue of the isomorphism $\cX \wedge \cO \cong \XoTDS$. Precisely, in the same manner that we had $\uF|n \ket = F_n | n \ket$, which led to the representation $\uFst = \sum F_n | n \ket$, presently we have $\utau | n \ket = \utau_n |n \ket$, so we have a state-like representation of $\utau$ given by $\utaust = \sum \utau_n |n \ket: \sum F_n |n \ket \Rarr \sum G_n |n \ket$. Projecting to any state $|n\ket$ we have $\bra n| \utau \ket : F_n = \bra n | \uF \ket \Rarr \bra n | \uG \ket = G_n$ so that $\bra n | \utau \ket = \utau_n$. Already at this point it is clear that if $\utau$ is a natural transformation, so is $\utau_0: F \Rarr G$, but we nevertheless work this out in full for illustrative purposes. Using the fact that $\utau$ is a natural transformation, since objects of $\cX$ are regarded as functors, for $f:A \rarr B$ in the base category of $\cX$, we have a commutative square in $\SetD \wedge \cO$ if objects of $\cX$ are valued in $\SetD$:
\beq
\xymatrix{
	\uF A \ar[d]_{\uF f} \ar[r]^{\utau_A} & \uG A \ar[d]^{\uG f} \\
	\uF B \ar[r]_{\utau_B} & \uG B
} \nonumber
\eeq
If $\uFst = \sum F_n | n \ket$, we have $|\uF \ket A = \sum F_nA |n \ket$ in $\SetD \oT^{\Delta} \Sigma$, and likewise $ \utaust_A = \sum (\utau_n)_A |n \ket$. 
Now because we have an isomorphism:
\begin{align}
	\cX \wedge \cO &\xrarr{\cong} \XoTDS \nonumber \\
	\uF & \mapsto \uFst \nonumber
\end{align}
being defined pointwise at objects $A$ of the base category of $\cX$, we have pointwise isomorphisms:
\begin{align}
	\SetD \wedge \cO &\xrarr{\cong} \SetD \oT^{\Delta} \Sigma \nonumber \\
	\uF A & \mapsto \uFst A \nonumber
\end{align}
The state-like representation of the above CD in $\SetD \wedge \cO$ expressed in $\SetD \oT^{\Delta} \Sigma$ using this isomorphism now reads:
\beq
\xymatrix{
	|\uF \ket A \ar[d]_{\uFst f } \ar[r]^{ \utaust_A} & | \uG \ket A \ar[d]^{|\uG \ket f } \\
	|\uF \ket B \ar[r]_{\utaust_B} & |\uG \ket B
} \nonumber
\eeq
which upon projection to a state $|0\ket$ produces:  
\beq
\xymatrix{
	\bra 0 | \uF \ket A \ar[d]_{\bra 0 | \uF \ket f} \ar[r]^{\bra 0 | \utau \ket_A} & \bra 0 | \uG \ket A \ar[d]^{\bra 0 | \uG \ket f} \\
	\bra 0 | \uF \ket B \ar[r]_{\bra 0 | \utau \ket_B} & \bra 0 | \uG \ket B
} \nonumber
\eeq
and with $F = \bra 0 | \uF \ket$, $Ff = \bra 0 | \uF \ket f$, $\tau = \bra 0 | \utau \ket$ this reads:
\beq
\xymatrix{
	FA \ar[d]_{Ff} \ar[r]^{\tau_A} & GA \ar[d]^{Gf} \\
	FB \ar[r]_{\tau_B} & GB
} \nonumber
\eeq
hence $\tau: F \Rarr G$.

\subsubsection{Maps on physical states}
We argued above that one should focus on those functors $\psi$ for which the map $F \rarr \psi F$ is effectively just changing the state the essential intrinsic functor $\uF$ initially finds itself in to another state. Thus we seek a map of physical states since this is primarily what we are interested in. We consider only those changes of states that are functors $\upsi: \Sigma \rarr \Sigma$. Let $\upsi$ be one such functor. We consider the following defining CD for $\psi$:
\beq
\xymatrix{
	\Sigma \oT (\cX \wedge \cO) \ar[d]_{\cong} \ar[r]^{\upsi \oT id} & \Sigma \oT (\cX \wedge \cO) \ar[d]^{\cong} \\
	\Sigma \oT (\XoTDS) \ar[d]_{m} \ar[r]^{\upsi \oT id} & \Sigma \oT (\XoTDS) \ar[d]^{m} \\
	\cX \ar[r]_{\psi} & \cX
}\nonumber
\eeq
with $\psi = m( \upsi \oT id)$. To fix ideas let $\upsi: |0 \ket \rarr |1 \ket$ in $\End (\Sigma)$, and take $\uF$ such that $\uF |0 \ket = F|0\ket$. Then we are looking at:
\beq
\xymatrix{
	\bra 0 | \oT \uF \ar[d]_{\cong} \ar[r]^{\upsi \oT id_{\uF}} & \bra 1 | \oT \uF \ar[d]^{\cong} \\
	\bra 0 | \oT \uFst \ar[d]_{m} \ar[r]^{\upsi \oT id_{\uFst}} & \bra 1 | \oT | \uF \ket \ar[d]^{m} \\
	\bra 0 | \uF \ket = F \ar[r]_-{\psi} & \bra 1 | \uF \ket = \bra 1 | \sum F_n |n \ket = \bra 1 | F_1 |1 \ket = F_1 = \psi F
}\nonumber
\eeq
This makes manifest the fact that $\psi$ is not strictly speaking a state in itself, rather it is derived from a functor of states $\upsi$. However the bottom horizontal map of the above diagram shows that $\psi$ shifts objects of $\cX$ from one state to another. In other terms we can regard $\psi$ as a relative state of our system, which explains why we nevertheless refer to $\dX$ as the Segal category of states for our system $\cX$.\\

That $\psi$ is a functor is due in part to the fact that $\upsi$ itself is a functor on $\Sigma$. Consider:
\beq
\xymatrix{
	\bra p | \ar[d]_{\uu} \ar[r]^-{\upsi} & \bra p' | \ar[d]^{\upsi \uu} \\
	\bra q | \ar[r]_-{\upsi} &\bra q' | 
}\label{CDpsifctr}
\eeq
evaluated at some $\uF \cong \uFst$ this yields:
\beq
\xymatrixcolsep{4pc}
\xymatrix{
	\bra p | \uF \ket \ar[d]_{m(\uu \oT id_{\uFst}) } \ar[r]^-{m(\upsi \oT id_{\uFst})} & \bra p' | \uF \ket \ar[d]^{m(\upsi \uu \oT id_{\uFst})} \\
	\bra q | \uF \ket  \ar[r]_-{m(\upsi  \oT id_{\uFst}) } &\bra q' | \uF \ket
}\label{CDpsifctrF}
\eeq
which reads:
\beq
\xymatrix{
	F_p \ar[d]_{u} \ar[r]^{\psi} & \psi F_p \ar[d]^{\psi u}\\
	F_q \ar[r]_{\psi} & \psi F_q
}\label{psifctrCDoneoftwo}
\eeq
with $u = m(\uu \oT id_{\uFst})$ and $\psi u = m( \upsi \uu \oT id_{\uFst})$. Now note that the notation $F = \bra 0 | \uF \ket$ is bifunctorial. A natural transformation $F \Rarr G$ in $\cX$ can originate from a change of states, as we just considered, or from a change in the underlying functor $\uF$ of $F$, something we call a \textbf{mutation}, $\utau: \uF \Rarr \uG$. Thus we also consider the following CD:
\beq
\xymatrixcolsep{4pc}
\xymatrix{
	\bra 0 | \oT \uFst \ar[d]_{id_{\bra 0 |} \oT \utaust = \bra 0| \oT \utaust} \ar[r]^-{ \upsi \oT id_{\uFst}} & \bra 1 | \oT \uFst \ar[d]^{id_{\bra 1 |} \oT \utaust = \bra 1 |  \oT \utaust} \\
	\bra 0 | \oT \uGst \ar[r]_-{\upsi \oT id_{\uGst}} & \bra 1 | \oT |\uG \ket
} \label{CDpsitau}
\eeq
which upon multiplication yields:
\beq
\xymatrixcolsep{5pc}
\xymatrix{
	F = \bra 0 | \uF \ket \ar[d]_{\tau = \bra 0 | \utau \ket} \ar[r]^{\psi = m (\upsi \oT id_{\uFst})} & \bra 1 | \uF \ket = \psi F \ar[d]^{\bra 1 | \utau \ket = \psi \tau} \\
G = \bra 0 | \uG \ket \ar[r]_{\psi = m(\upsi \oT id_{\uGst})} & \bra 1 | \uG \ket = \psi G
}\label{CDpsitaum}
\eeq
which simplifies as:
\beq
\xymatrix{
	F \ar[d]_{\tau} \ar[r]^{\psi} & \psi F \ar[d]^{\psi \tau} \\
	G \ar[r]_{\psi} & \psi G
} \label{psifctrCDtwooftwo}
\eeq
hence \eqref{psifctrCDoneoftwo} and \eqref{psifctrCDtwooftwo} show that $\psi$ is a functor in $\RHomXX$.\\

We pointed out that $F = \bra 0| \uF \ket$ is bifunctorial. Note that this really follows from $m$ being bifunctorial on $\Sigma \oT (\XoTDS)$ as $F = \bra 0 | \uF \ket = m(\bra 0 | \oT \uFst)$. Using this line of argument, right tensoring \eqref{CDpsifctr} by $\uFst$ in $\Sigma \oT (\XoTDS)$ and using the fact that $m$ is functorial in its first argument, one obtains \eqref{CDpsifctrF}. Then using the bifunctoriality of $m$ again, one goes from \eqref{CDpsitau} to \eqref{CDpsitaum}. Thus $\psi$ being a functor follows from the bifunctoriality of $m$, and the fact that $\upsi$ is a functor.

\subsection{Iterative shifts}
So far applying $\psi$ to some $F \in \cX$ just produces a shift in $F$. Dynamics comes up when one considers a sequence of such shifts. For simplicity, we initially focus on dynamics generated by a single $\psi \in \dX$. Consider a sequence of states in $\Sigma$ generated by $\upsi$ from a physical state $|0 \ket$ as in:
\beq
| 0 \ket \xrarr{\upsi} \upsi |0\ket = | 1 \ket \xrarr{\upsi} \upsi |1 \ket = \upsi^2 |0 \ket \xrarr{\upsi} \cdots \nonumber
\eeq
which once tensored with $\uFst$ and then evaluated produces:
\beq
\xymatrixcolsep{4pc}
\xymatrix{
	\bra 0 | \uF \ket \ar@{=}[ddd] \ar[r]^{m(\upsi \oT id_{\uFst})} & \bra 1 | \uF \ket  \ar@{=}[d] \ar[r]^{m(\upsi \oT id_{ \uFst})} & \bra 2 | \uF \ket \ar@{=}[d] \ar[r] & \cdots \\
&\bra 0 | \upsi | \uF \ket \ar@{=}[dd] &   \bra 1| \upsi | \uF \ket \ar@{=}[d] \\
&& \bra 0 | \upsi^2 | \uF \ket \ar@{=}[d]\\ 
	F \ar[r]_{\psi} & \psi F \ar[r]_{\psi} &\psi^2 F \ar[r]_{\psi} &\cdots
} \nonumber
\eeq

As argued before, the next step consists in considering states that vary via a functor $\Psi: \dX \rarr \dX$, thereby producing a sequence in $\dX$, and take the colimit thereof. Let us denote by $\eps_0$ an initial state in $\dX$. We consider:
\beq
\xymatrix{
	\eps_0 \ar[r]^-{\Psi} &\Psi \eps_0  \ar[drrr]_{\iota_1} \ar[r]^-{\Psi} & \Psi  \eps_1  =  \Psi^2 \eps_0  \ar[drr]^{\iota_2} \ar[r] & \cdots \ar@{.>}[dr]  \\
	&&&&  \colim \Psi^n \eps_0 
} \nonumber
\eeq
which exists since $\cX$ being a Segal topos, so is $\dX = \RHomXX$, hence it has all small colimits. One defines a \textbf{flow} of our system from an initial state $\eps_0$ and induced by some $\Psi \in \RHom(\dX, \dX)$ to be given by $\Psi^{\infty} \eps_0  = \colim_n \Psi^n \eps_0 $. This is defined pointwise. Thus for some object $F \in \cX$, a flow in the dynamics of $F$ starting from an initial state $\eps_0$ in the direction of $\Psi$ is given by an augmented diagram:
\beq
\xymatrix{
	F \ar[r]^{\eps_0} & \eps_0F \ar[r]^-{\Psi} & \Psi \eps_0 F   \ar[drrr] \ar[r]^{\Psi} & \Psi^2  \eps_0F \ar[drr]   \ar[r] &  \cdots \ar@{.>}[dr]  \\
	&&&&& \Psi^{\infty} \eps_0F   
} \nonumber
\eeq

To summarize the current picture, our system $\cX$ is identified with $id_{\cX} \in \RHomXX:= \dX$, a Segal category of states of our system, which is enhanced to $\RHom(\dX, \dX):= \dtwoX$, and objects thereof are the functors $\Psi$ discussed above, which effectively move our system from one state to another. Iterative compositions of a given $\Psi$ on an initial state ultimately produce a flow of our system in the limit. The repetitiveness of using $\psi \in \dX$ above was fixed by introducing $\Psi \in \dtwoX$. Now however that problem is repeating itself for $\Psi$. Thus iterating this overall construct, and letting $\RHom(\dnX, \dnX):= \delta^{n+1} \cX$ represent higher states, we then claim the $\infty$-Segal topos $\colim_{n \infty} \dnX$ - colimit taken in the 2-Segal category of Segal categories - contains all possible higher dynamics of our initial system $\cX$, and consequently corresponds to a complete description thereof, making $\cX$ a closed system in the physical sense of the word. To fix ideas, the first stages of this construction go as follows. If $\psi^{(-,m)} \in \delta^m \cX$ and $\psi^{(n+1,m)} = \psi^{(n,m+1)} \psinm$, we get a discrete flow at the level of $\cX$:
\beq
F \xrarr{\psi^{(0,1)}} \psi^{(0,1)} F \xrarr{\psi^{(1,1)}} \psi^{(1,1)} \psi^{(0,1)}F \xrarr{\psi^{(2,1)}} \psi^{(2,1)} \psi^{(1,1)} \psi^{(0,1)} F \rarr \cdots \nonumber
\eeq
where the objects $\psi^{(n,1)}$ of $\dX$ fit into a sequence:
\beq
\psi^{(0,1)} \xrarr{\psi^{(0,2)}} \psi^{(0,2)} \psi^{(0,1)} = \psi^{(1,1)} \xrarr{\psi^{(1,2)}} \psi^{(1,2)} \psi^{(1,1)} = \psi^{(2,1)} \rarr \cdots \nonumber
\eeq
and the objects $\psi^{(n,2)}$ of $\dtwoX$ themselves fit into a sequence:
\beq
\psi^{(0,2)} \xrarr{\psi^{(0,3)}} \psi^{(0,3)} \psi^{(0,2)} = \psi^{(1,2)} \xrarr{\psi^{(1,3)}} \psi^{(1,3)} \psi^{(1,2)} = \psi^{(2,2)} \rarr \cdots \nonumber
\eeq
and higher states are constructed in like manner.\\

\subsection{Remarks on higher states}
To come back to the iterative step $\delta^{n+1} \cX = \RuHom(\dnX, \dnX)$, an object of $\delta^{n+1} \cX$ is a functor of Segal categories $\dnX \rarr \dnX$ that provides a global shift in the objects of $\dnX$. The point we want to make is that $\Psi_{n+1} \in (\delta^{n+1} \cX)_0$ is but an object of $\delta^{n+1} \cX$, yet from the perspective of $\dnX$, for all $\zeta \in (\dnX)_0$, $\Psi_{n+1}: \zeta \mapsto \Psi_{n+1}\zeta$, in other terms its action can be quite complex. Essentially then simple objects in a higher category can have complex realizations in a lower category, or viewed differently, complex operations in a given categorical context can be repackaged into something far simpler in a higher category. This is especially pertinent in physical applications when faced with local behaviors at the level of $\dnX$, which may point to the existence of a possible overarching theory without providing specifics about such a theory for lack of clarity. Expanding our views at the level of $\dnX$ to something global can provide the object of $\delta^{n+1} \cX$, if one can extract it, that explains the overall picture, and local behaviors in particular.\\

Another interesting observation is the following: if we place ourselves at the level of $\dnX$ and consider $\zeta$ an object of $\dnX$, and $\Psi \in \RuHom(\dnX, \dnX) = \delta^{n+1}\cX$ is a higher state of our system, from the perspective of $\dnX$, $\zeta$ and $\Psi \zeta$ correspond to a ``before" and an ``after", hence a de facto notion of time, while at the level of $\delta^{n+1}\cX$ already, objects $\Psi$ thereof are defined pointwise, so $\Psi$ is defined by $\zeta \xrarr{\Psi} \Psi \zeta$, hence regard past (or present) and future events in $\dnX$ as being simultaneous.\\

\subsection{Quantum states} \label{qs}
Note that we have defined a generic flow as being defined from higher states $\psinm \in \dmX $ for $m \geq 1$, $n \geq 0$. In practice though one would have to consider only those initial states that are fundamental, such as quantum states. In the present section, we define such fundamental, quantum states in our formalism. Since this is fairly applied, we will work with $\cX = \dAffCtildetau$ for definiteness. Let $\psi \in \dX$. For $ \psi$ to be a non-trivial quantum state, we first ask that there is some $F \in \cX$, $\psi F \neq F$, that is there is some $A \in \CommC$, $\psi FA \neq FA$. We define the \textbf{differential matrix} of $\psi$ to be defined by:
\beq
\Delta \psi = \big(d_{\SetD}(\psi FA,FA)\big)_{\substack{\!\!\!F \in \cX \\  A \in \CommC}} \nonumber
\eeq
and for simplicity we can write $(\Delta \psi)_{F,A} = d_{\SetD}(\psi FA, FA)$. Using the fact that the geometric realization of a simplicial set is a CW-complex, for $X,Y \in \SetD$, we define:
\beq
d_{\SetD}(X,Y) = d_{CW}(|X|,|Y|) \nonumber
\eeq
where we consider homotopy types of CW-complexes, which are determined by their homotopy classes of attaching maps, thus we can define a distance function $d_{CW}$ on CW-complexes, denoted $d_{CW} (C,C')$ for $C$ and $C'$ two CW-complexes, by being the number of differing homotopy classes of attaching maps $C$ and $C'$ have. This number may be unbounded, but this is inconsequential for the sake of defining quantum states, whose differential matrices have finite entries. Denote by $\delta(C_n,C_n')$ the number of differing homotopy classes of attaching maps on the n-skeleta $C_n$ and $C_n'$ of $C$ and $C'$ respectively to get $C_{n+1}$ and $C_{n+1}'$, and let:
\beq
\delta(C_{\leq n},C_{\leq n}') = \sum_{0 \leq p \leq n} \delta(C_p,C_p') \nonumber
\eeq
where $C_{\leq n} = \bigcup_{0 \leq p \leq n} C_p$. Using this, we let:
\beq
d_{CW}(C,C') = \sum_{n \geq 0}  \delta(C_n,C_n') = \lim_{n \infty} \delta( C_{\leq n}, C'_{\leq n}) \nonumber
\eeq
This is clearly symmetric, reflexive, and the only situation where we would run into difficulties while showing $\delta(C_n,C_n'') \leq \delta(C_n,C_n') + \delta(C_n',C_n'')$ is when the left hand side is 1, and the right hand side is $0 + 0$, which cannot occur by virtue of the fact that homotopy is an equivalence relation, hence transitive. It follows that $d_{CW}$ is a well-defined distance on CW-complexes. Having defined this, it follows that $d_{\SetD}(X,Y) = d_{CW}(|X|,|Y|)$ is a well-defined distance function on simplicial sets in a homotopical context. Hence $(\Delta \psi)_{F,A} = d_{\SetD}(\psi FA,FA)$ is well-defined. We now define \textbf{quantum states} to be those states $\psi$ that make their differential matrix minimal. To be precise, $\psi$ is a quantum state if there is no other state $\eta$ such that:
\beq
(\Delta \eta)_{F,A} < (\Delta \psi)_{F,A} \nonumber
\eeq
for some $F \in \cX$, and some $A \in \CommC$, all other entries of $\Delta \eta$ being higher than, or equal to those of $\Delta \psi$ otherwise.\\

We can iterate this for higher states. Consider $\Psi \in \dtwoX$. It is a nontrivial (quantum) state if there exists $\eta \in \dX$, $\Psi \eta \neq \eta$, which holds if there exists $F \in \cX$, $\Psi \eta F \neq \eta F$, and again for objects $F$ that are functors, this is true if there exists $A \in \CommC$ such that $\Psi \eta F A \neq  \eta F A $ in $\SetD$. We are led to considering the differential matrix of $\Psi$:
\beq
\Delta \Psi = \big( d_{\SetD}(\Psi \eta F A, \eta F A) \big)_{\eta, F, A} \nonumber
\eeq
and states $\Psi \in \dtwoX$ that are quantum in nature are those that make this cubic matrix minimal. More generally, physical, smooth flows are built from higher states $\psinm \in \delta^m \cX$ that must all be quantum in nature, each having a minimal differential matrix of size $m+1$.

\section{Cusps} \label{Cusps}
A natural question is, given a state $\eps $ in $\dX$, and $\Psi: \dX \rarr \dX$, what are the other such functors $\Xi$ that yield a same local flow as $\Psi$ immediately at $\eps $. We claim those are satisfying $T_{\Xi,\eps} \cong T_{\Psi,\eps}$ for a same state $\eps$ in $\dX$, an isomorphism of tangent spaces. Since one of the things we want to achieve is being able to clarify the functorial dynamics at the level of a Segal category $\cX$ such as $\dAffCtildetau$, we focus on this particular example presently, and we take for our base monoidal category a symmetric monoidal model category $\cC$ that we choose to be simplicial, that is we place ourselves within the realm of derived algebraic geometry. The correct notion of tangent space in derived algebraic geometry being that of a tangent complex, we define $T_{\Psi,\eps} := \beeT_{\Psi,\eps}$ as the proper adaptation in our context of the tangent complex of $\Psi$ at $\eps $. Then we define the \textbf{cusp} of $\Psi$ at $\eps $ by:
\beq
\Cusp[\Psi,\eps] = \{ \Xi: \dX \rarr \dX \, | \, \beeT_{\Xi,\eps} \cong \beeT_{\Psi,\eps}  \} \nonumber
\eeq
it gives those functors $\Xi: \dX \rarr \dX$ that have a same flow as $\Psi$ at $\eps $. Having thus defined cusps, note that this defines equivalence classes of local flows in the space of states; from any given state $\eps $, functors $\Psi: \dX \rarr \dX$ can be grouped in equivalence classes $[ \Psi ] := \Cusp[\Psi, \eps]$.

\subsection{Formal construction of the tangent complex}
We remind the reader of the formal construction of a tangent complex as covered in \cite{TV4}. We will not need this construction in our work, but it is nevertheless important to understand how the tangent complex comes about to justify our present definition of the tangent object $\TPsix$. For $F \in \AffCtildetau$ a stack, $x: \RuSpec A \rarr F$ an $A$-point of $F$, $\TFx$ is an element of the homotopy category of spectra $\Ho(\Sp(\AMod))$. Here $A \in \cA$, $\cA \subset \CommC$ having $\cC_0$-good objects for elements, where according to Def 1.1.0.10 of \cite{TV4}, $A \in \CommC$ is good with respect to $\cC_0$ if the functor $\Ho(\AMod^{\op}) \rarr \Ho((\AMod_0^{\op})^{\,\wedge})$ is fully faithful, where $\AMod_0 = \AMod \cap \cC_0$ and $\cC_0$ is a full subcategory of $\cC$ satisfying Assumption 1.1.0.6 of \cite{TV4} (see reference for more details on why it is needed). Here we have used the notation $\RuSpec A = \bR \underline{h}_{\Spec A}$, for $\bR \underline{h}: \Ho(\AffC) \rarr \St(\cC, \tau)$, the right derived functor of the model Yoneda embedding. In this definition, for $\Gamma_*$ a cofibrant resolution functor on the model category $\AffC$, $\underline{h}_X(Y) = \Hom_{\AffC}(\Gamma_*(Y),X)$.\\

With regards to spectra, one uses the fact that $\cC$ being pointed, one has a suspension functor $S: \Ho(\cC) \rarr \Ho(\cC), x \mapsto S(x) = * \coprod_x^{\bL} *$. Fix once and for all, as done in \cite{TV4}, an object $S^1_{\cC}$ in $\cC$, a cofibrant model for $S(1) \in \Ho(\cC)$, write $S^1_A = S^1_{\cC} \oT A \in \AMod$. That allows one to consider transition morphisms $\sigma_n: S^1_A \oTA M_n \rarr M_{n+1}$. Following the same argument as in \cite{TV4}, one gets a category of spectra $\Sp(\AMod)$. Put the stable model structure on $\Sp(\AMod)$ by taking the left Bousfield localization of this category, using as local objects stable modules $M_* \in \Sp(\AMod)$, objects such that the induced morphisms $M_n \rarr \RHom_A(S^1_A, M_{n+1})$ are isomorphisms in $\Ho(\AMod)$, where $\uHom_A(S^1_A,-)$ is the right adjoint of $S^1_A \oTA -$. Further $\Ho(\Sp(\AMod))$ inherits from $\Ho(\AMod)$ a symmetric monoidal structure $-\oTAL-$ making it into a closed symmetric monoidal category with an internal hom $\RHom_A^{\Sp}(M_*, N_*) \in \Ho(\Sp(\AMod))$ for $M_*$ and $N_*$ two stable $A$-modules. This allows one to define:
\beq
\TFx = \RHom_A^{\Sp}(\LFx,A) \nonumber
\eeq
for $\LFx \in \Ho(\Sp(\AMod))$ the cotangent complex of $F$ at $x$. The latter is defined by:
\beq
\Der_F(X,-) \cong \bR \uh_s^{\LFx} \in \Ho((\AMod_0^{\op})^{\, \wedge}) \nonumber
\eeq
with $\Der_F(X,M) = \Map(X[M],F) \in \Ho(\SetD)$, $X[M] = \RuSpec (A \oP M)$, $A \oP M$ the square zero extension of $A$ by $M \in \AMod$. This defining isomorphism also makes use of a functor $\uh_s: \Sp(\AMod)^{\op} \rarr (\AMod_0^{\op})^{\, \wedge}$, with $\uh_s^{M_*}: \AMod_0 \rarr \SetD$ mapping $N$ to $\Hom(M_*, \Gamma_*(S_A(N)))$, $\Gamma_*$ a simplicial resolution functor on $\Sp(\AMod)$, and the functor $S_A: \Ho(\AMod) \rarr \Ho(\Sp(\AMod))$ is defined by $S_A(M)_n = (S^1_A)^{\oTA^n} \oTA M$. One has that for all $A \in \CommC$, $\bR \uh_s: \Ho(\Sp(\AMod))^{\op} \rarr \Ho((\AMod_0^{\op})^{\, \wedge})$ is well-defined as the right derived functor of $\uh_s$. All of this work can be found in \cite{TV4}.\\

Perhaps we can pause at this point and note a first departure on our end from this standard formalism. Stacks are valued in simplicial sets, but our functors, to which we apply this construction, are valued in Segal categories. Indeed, $\Psi: \dX \rarr \dX$ for us plays the role of $F$ above. Thus the definition of $\Der_F(X,-)$ will have to be modified in our situation (see \cite{TV4} for the discussion about this functor having a lift to a functor $\AMod \rarr \SetD$). Additionally, $\Map$  is the mapping complex for a model category, which here in particular is $\AffCtildetau$, and in our case we work with $\RHom(\dX, \dX)$, for which we don't have an obvious mapping complex. We will remedy this below.\\

\subsection{Derivations}
Letting $\dtwoX = \RHom(\dX, \dX) \in \Ho(\SeCat)$, and supposing for a moment for $A \in \dX$ we have a functor $\RuSpec A \in \dtwoX$, then a natural transformation $x: X=\RuSpec A \rarr \Psi$ in $\dtwoX$ providing an $A$-point of $\Psi$ would be an element of $(\dtwoX)_{(\RuSpec A, \Psi)} \in \Ho(\SetD)$. However our base category $\dX$ is no symmetric monoidal category, so we have no obvious notion of linear algebra over it, thus a priori one cannot talk about $A$-modules, let alone square zero extensions $A \oP M$. Note however that what we are looking for is deformations of $A$, and this can be implemented with infinitesimal compositions $MA$ for some functor $M \in \dX$. We therefore define:
\beq
\Der_{\Psi}(X,M) = \Map(\RuSpec(MA), \Psi) \nonumber
\eeq
We first have to formally define $\RuSpec A$, which will lead us to defining certain mapping objects in generalized categories. Having done that, the definition of the analog in our context of the cotangent and tangent complexes will fall out of our formalism.\\

\subsection{Mapping spaces}
Presently we discuss $\RuSpec A \in \dtwoX$ for $A \in \dX$, and what that entails for morphisms of $\dtwoX$ that have it for domain. For $B \in \dX$, we define:
\beq
\RuSpec A(B) = \Nat(A, B) \nonumber
\eeq
to be the collection of natural transformations $\Lambda_{i,B}: A \Rarr B$. To make sense of a morphism $\Omega_B: \RuSpec A(B) \rarr \Psi(B)$, we therefore have to look at $\Psi$ differently. $\Psi(B) \in \dX = \RHomXX$ being a functor from $\cX$ to itself, for $f: F \rarr G$ in $\cX$, we have a commutative diagram:
\beq
\xymatrix{
	F \ar[d]_f \ar[r] & \Psi(B)F \ar[d]^{\Psi(B)f} \\
	G \ar[r] &\Psi(B)G
} \nonumber
\eeq
which we can regard as defining a natural transformation $\Gamma_B: id \Rarr \Psi(B)$:
\beq
\xymatrix{
	id(F) \ar[d]_{id(f)} \ar[r]^{\Gamma_{B,F}} & \Psi(B)F \ar[d]^{\Psi(B)f} \\
	id(G) \ar[r]_{\Gamma_{B,G}} &\Psi(B)G
} \nonumber
\eeq
in a natural way, or equivalently $\Psi: \dX \rarr \dX$ can be viewed as a natural transformation valued functor $\underline{id} \Rarr \Psi$, where $\underline{id}$ is identically the identity functor on $\dX$. In this manner, a morphism $\Lambda_{i,B} \rarr \Psi(B)$ is regarded as a 2-morphism $\Omega_{i,B}: \Lambda_{i,B} \rarr \Gamma_B$ in $\dX$ with $(\Lambda_{i,B}: A \Rarr B) \rarr (\Gamma_B: id \Rarr \Psi(B))$ as presentation, and the collection of all such maps gives us $\{ \Nat(A, B) = \RuSpec A(B) \rarr (id \Rarr \Psi(B)) \}$. We define $\Map( \RuSpec A, \Psi)$ pointwise in this manner. Observe that this gives us a notion of an A-multi-point of $\Psi$.\\

\subsection{Generalized categories and generalized functors}
One pressing issue at this point is making sense of the collection of morphisms $\Der_{\Psi}(X,M)(B) = \Map(\RuSpec MA(B),\Psi(B))$ mentioned above. We start by defining generalized categories, and then we define the corresponding notion of generalized functors.\\

In an ordinary category $\cC$, an object $c$ of $\cC$ is an element of the set $\Ob(\cC)$ of objects of $\cC$. We generalize this to considering the case where we have \textbf{amalgams} of objects $c = (c_1, \cdots, c_n ) \in \Ob(\cC)^n$. We are led to introducing the generalized category $\cC^{gen}$ induced by $\cC$ as the graded object $\cC^{gen} = \oP_{n \geq 1} \cC^{gen}_n$ with $\Ob(\cC^{gen}_n) \cong (\Ob \cC)^n$. The appropriate notion of morphism in this case is defined as follows: given objects $c = (c_1, \cdots, c_n ) \in \cC^{gen}_n$ and $d = (d_1, \cdots, d_m ) \in \cC^{gen}_m$, a morphism $f$ from $c$ to $d$ is given by a matrix of morphisms: $f = (f_{ij}: c_i \rarr d_j ) \in \Mor(\cC)^{n \times m}$. More generally, denoting such a morphism by $f_{nm}$, a generic morphism in $\cC^{gen} = \oP \cC^{gen}_n$ is of the form $f = \oP f_{nm}$. Coming back to the special case where $f = (f_{ij}):c \rarr d$ above, to account for the fact that for $i$ and $j$ fixed, there may be several morphism between $c_i$ and $d_j$, we denote such morphisms by $f_{ij}^{\lambda}:c_i \rarr d_j$, with $\lambda \in A_{ij}:=A[c_i, d_j]$ an indexing set for such morphisms, whose elements are referred to as \textbf{flavors}. In the present formalism we regard $f_{ij}^{A_{ij}} = (f_{ij}^{\lambda})_{\lambda \in A_{ij}}$ as a \textbf{quiver} (\cite{PG}) between $c_i$ and $d_j$. Write:
\beq
f^A  = (f_{ij}^{A_{ij}}) = ((f_{ij}^{\lambda})_{\lambda \in A_{ij}})_{\substack{1 \leq i \leq n \\ 1 \leq j \leq m}} := (f_{ij}^{\lambda})  \in \Mor(\cC^{gen}_n, \cC^{gen}_m) \nonumber
\eeq
with $A = \times A_{ij}$, and we refer to $f^A$ as a (higher-dimensional) \textbf{cube} in reference to the fact that we have three indices $i$, $j$ and $\lambda$ (and later as we will soon see, multi-flavors). By abuse of notation, we will just write $f = (f^{\lambda}_{ij})$, $A$ being implied. In the immediate, for simplicity, we focus on morphisms that are not made of quivers. For $e = (e_1, \cdots, e_p )$ another object, and $g:d \rarr e$ being given by $g = ( g_{jk}: d_j \rarr e_k )$, the composition of $g$ and $f$ is defined by $g \circ f = (g_{jk} \circ f_{ij}: c_i \rarr e_k )$. It is clearly associative, with unit $id_c = ( id_{c_i} )$. It is noteworthy that with this definition $(g \circ f)_{ik}:c_i \rarr e_k$ has a flavor $j$, keeping track of what object $d_j$ is $(g \circ f)_{ik}$ going through, and we write $(g \circ f)_{ik}^j$ for such a morphism. We now add flavors to the picture: let $g = (g_{jk}^{\mu}: d_j \rarr e_k)$ and $f = (f_{ij}^{\lambda}: c_i \rarr d_j)$. The composite $g \circ f$ will now have flavors $j$, and for $j$ fixed, flavors $\lambda \in A[c_i,d_j]$ and $\mu \in A[d_j,e_k]$. For $j$ fixed, we therefore have two flavors $\lambda$ and $\mu$. If one is fixed, the other one acts as a flavor relative to the first one, as follows:
\beq
(g \circ f)_{ik}^{j, \lambda} = ((g \circ f)_{ik}^{j, \lambda, \mu})_{\mu \in A[d_j,e_k]}:= ((g \circ f)_{ik}^{j, \lambda, \mu}) \nonumber
\eeq
so that:
\begin{align}
	g \circ f &= ((g \circ f)_{ik}) \nonumber \\
	&= ( ( (g \circ f)_{ik}^j)  )  \nonumber \\
	&= ( \cdots ( g \circ f)_{ik}^{j, \lambda} \cdots ) \nonumber \\
	&= ( \cdots (g \circ f)_{ik}^{j,\lambda, \mu} \cdots ) \nonumber \\
	&= ( (g \circ f)_{ik}^{j,\lambda, \mu}) \nonumber
\end{align}
where at each level the new index is a flavor relative to the indices of the one line prior. One could have fixed $\mu$ and obtained $\lambda$ as a flavor relative to $\mu$, for a same end result entry-wise. The fact that we have towers of flavors justifies calling such objects higher dimensional cubes. Note however that with this notion of generalized category, one cannot define inverse morphisms without running into pathological behavior. More generally still, a \textbf{generalized category} is a graded object $\cC = \oP \cC_n$ with $\cC_n \cong \prod_{1 \leq i \leq n} \cC_{(i)}$, the $\cC_{(i)}$'s being ordinary categories, and morphisms in $\cC$ are defined levelwise from objects $c = (c_1,\cdots, c_n) \in \cC_n$ to $d = (d_1,\cdots,d_m) \in \cC_m$ by $f = (f_{ij})$ with $f_{ij} \in \Hom(\cC_{(i)},\cC_{(j)})$. Finally note that if $\cC$ is a generalized category, so is $\Mor \cC$, hence we can apply this formalism to define higher morphisms in a generalized setting.\\

For the sake of completeness, we define the notion of \textbf{generalized functor}, which are functors between generalized categories. The simplest case is the following: any ordinary functor $F: \cC \rarr \cD$ between ordinary categories gives rise in a natural way to a generalized functor $F^{gen}: \cC^{gen} \rarr \cD^{gen}$ defined by $F^{gen} = \oP_{n \geq 1} F^{gen}_n$, with $F^{gen}_n: \cC_n \rarr \cD_n$. By abuse of notation we will just denote $F^{gen}_n$ by $F$. Then if $c = (c_1, \cdots, c_n) \in \cC^{gen}_n$, $Fc = (Fc_1, \cdots, Fc_n )$. If $d = (d_1, \cdots, d_m )$ is another object of $\cCgen$, $f: c \rarr d$ is a morphism given by $f = ( f_{ij}: c_i \rarr d_j )$, then $Ff = (F(f_{ij}): Fc_i \rarr Fd_j )$, that is if $f \in \Hom(\cC^{gen}_n, \cC^{gen}_m)$, then $Ff \in \Hom(\cD^{gen}_n, \cD^{gen}_m)$. If $e = (e_1, \cdots, e_p )$ is another object of $\cCgen$, and $g: d \rarr e$ is a morphism given by $g = (g_{jk}: d_j \rarr e_k )$, then:
\begin{align}
	Fg \circ Ff &= (F(g_{jk}) \circ F(f_{ij}) ) \nonumber \\
	&= ( F(g_{jk} \circ f_{ij}): Fc_i \rarr Fe_k ) \nonumber \\
	&= ( F((g \circ f)_{ik}^j )) \nonumber\\
	&= F((g \circ f)^j_{ik}) \nonumber \\
	&=F(g \circ f) \nonumber
\end{align}
We also have:
\beq
F(id_c) = F((id_{c_i} )) = (F(id_{c_i}) ) = (id_{Fc_i} ) = id_{Fc} \nonumber
\eeq

Presently we have not used the fact that the entries of $f$ and $g$ could have been quivers, but adding in flavors is just a tedious exercise, yielding the same results as above. Functors between generic generalized categories can be defined in an obvious way from the above work.\\

An equivalent definition of generalized categories consists in using sets to represent amalgams and morphisms: $c = \{c_1, \cdots, c_n \}$, and $f: c \rarr d$ is given by $f = \{f_{ij}: c_i \rarr d_j \}$. If $c_i \in \Ob(\cC)$, $\cC$ an ordinary category, $c$ is an object of $\cCgen$, and with $f_{ij} \in \cC_1$, we have $f \in \cCgen_1$, but note that in both instances we are dealing with sets, while in the cubic picture, for $c = (c_1,\cdots, c_n)$ and $d = (d_1, \cdots, d_m)$, $f = (f_{ij}) \in \cCgen_1$ is an object of $(\cC_1)^{n \times m}$, and if we have multiple flavors $\lambda_{k_{p+1}} \in A_{ij,\lambda_{k_1}, \cdots, \lambda_{k_p}}$, then $f = (f_{ij}^{\lambda_{k_1}, \cdots, \lambda_{k_r}}) \in (\cC_1)^N$, where $N = \sum_{ij}|A_{ij}| \times \cdots \times |A_{ij,\lambda_{k_1}, \cdots, \lambda_{k_{r-1}}}|$. To make a distinction between concepts appearing in both the Set-based formalism of generalized categories and its cubic counterpart, we will use a $\sigma$ upperscript in the latter case. We have initially given a cubical approach to generalized categories since it is new, and we have given a short Set-based description since it is most natural. In what follows, we will use the Set-based approach for the same reason, and later we will present the cubical approach, which we refer to as the simplicial case for reasons that will then become obvious.\\

To come back to our problem of formally constructing $\Der_{\Psi}(X,M)$, in the Set-based setting one can write $\Lambda_B^M = \{ \Lambda_{i,B}^M: MA \Rarr B \}$, $i$ being regarded as a flavor here, thereby casting $\dX$ as a generalized Segal category, albeit with just quivers on single objects, and $\Lambda^M_B = \Nat(MA,B)$ as a generalized morphism in $\dXgen$. Observe that if $\Lambda_{i,B}^M \in \dX_{(MA,B)} \in \Ho(\SetD)$, in a generalized context $\Lambda_B^M \in (\dXgen)_{(MA,B)}$ is an object of $\Set$. We also have $\Gamma_B: id \Rarr \Psi(B)$ as another morphism in $\dX$. The various 2-morphisms between $\Lambda^M_B$ and $\Gamma_B$ are given entrywise by some $\Omega_{i,B}^{M, \lambda}: \Lambda_{i,B}^M \rarr \Gamma_B \in \dX_{(MA,id, \Psi(B))} \in \Ho(\SetD)$ - as we will show in the next subsection - for $\lambda$ a flavor in some indexing set $L_i$. A collection of such maps gives a (quiver) morphism $\Omega^{M,L_i}_{i,B} = \{ \Omega^{M, \lambda}_{i,B} \, | \, \lambda \in L_i \}$, and such morphisms themselves can be collected into a morphism $\Omega_B^{M, L} = \{ \Omega_{i,B}^{M, L_i}: \Lambda_{i,B}^M \rarr \Gamma_B \}$ with $L = \times_i L_i$, an object of $ (\dXgen)_{(MA,id, \Psi(B))}$. Note that this gives us a generalized morphism $\Omega_B^{M,L}: \Lambda_B^M \rarr \Gamma_B$ if each index $i$ is accounted for in such a collection. Henceforth we will only consider such \textbf{``full"} collections. Using the fact that:
\beq
 \Lambda_B^M = \{\Lambda_{i,B}^M : MA \rarr B \} = \Map(MA,B) = \RuSpec(MA)(B) \nonumber
\eeq
it follows that one can define:
\begin{align}
	\Omega_B^M &= \{ \Omega^{M,L}_B: \Lambda^M_B \rarr \Gamma_B \} \nonumber \\
	&=  \Map(\RuSpec(MA)(B), \Psi(B) ) \nonumber \\
	&= \Der_{\Psi}(X,M)(B)  \nonumber
\end{align}

\subsection{Cotangent set}
Following the above philosophy of regarding $\delta A = MA$ as a small deformation of $A \in \dX$, one can naively expect the analog in our situation of the cotangent complex of $\Psi$ at $A$ to be an object of the form $\Map(\Psi A , \Psi(MA))$. This will help us see that the abstract definition of the cotangent object $\LPsix$ as we will define it does what we want it to do. For now, observe that $\Psi: \dX \rarr \dX$ is a functor on $\dX$, whose objects are states, which are really functors. Being functors from $\cX$ to $\cX$, they can be composed. If $A$ and $B$ are two such functors, one can define their composite $BA$. For the sake of consistency, we ask that $\Psi$ be a morphism on $(\dX, \circ, id)$, so that $\Psi(BA) = \Psi(B) \circ \Psi(A)$. Thus we are looking at $\Nat(\Psi A, \Psi M \Psi A)$. We regard $\Psi M$ as providing an infinitesimal deformation of $\Psi A$ in the $M$ direction.\\

As a first step towards defining a cotangent object, we now justify $\Omega^{M,\lambda}_{i,B}: \Lambda^M_{i,B} \rarr \Gamma_B$ being valued in $ \dX_{(MA,id,\Psi(B))} \in \Ho(\SetD)$ as initially mentioned in the preceding section. To put things in perspective, recall that we have $\Der_{\Psi}(X,M)(B) = \Map(\RuSpec (MA)(B), \Psi(B))$ is a set of full morphisms $\Omega_B^{M, L }$ of the form $\{ \Omega_{i,B}^{M, L_i}: \Lambda_{i,B}^M \rarr \Gamma_B\}$, with $L = \times_i L_i$, and where $\Lambda_{i,B}^M: MA \Rarr B$ and $\Gamma_B: id \Rarr \Psi(B)$. The quivers $\Omega^{M,L_i}_{i,B}$ are collections of maps $\Omega_{i,B}^{M,\lambda}$ between natural transformations, given by commutative diagrams:
\beq
\xymatrix{
	MA \ar[d]_{\Lambda_{i,B}^M} \ar[r] & id \ar[d]^{\Gamma_B} \\
	B \ar[r] & \Psi(B)
} \nonumber
\eeq
which because they are commutative diagrams, can be represented as:
\beq
\xymatrix{
	MA \ar[d] \ar@{.>}[dr] \ar[r] & id \ar[d] \\
	B \ar[r] & \Psi(B)
} \nonumber
\eeq
that is $\Omega_{i,B}^{M,\lambda}: MA \rarr id \rarr \Psi(B)$, or equivalently $\Omega_{i,B}^{M, \lambda}: MA \rarr B \rarr \Psi(B)$, thus $\Omega_{i,B}^{M, \lambda} \in \dX_{(MA,id,\Psi B)}$ or $\Omega_{i,B}^{M, \lambda} \in \dX_{(MA,B,\Psi B)}$ depending on which factorization one is looking at. $\dX$ being a Segal category, we therefore have $\Omega_{i,B}^{M, \lambda} \in \Ho(\SetD)$. Note however that $\Omega_B^{M,L_i} = \{ \Omega_{i,B}^{M, \lambda}: \Lambda_{i,B}^M \rarr \Gamma_B \, | \, \lambda \in L_i \} \in  (\dX)^{gen}_{(MA,id,\Psi(B))}$ is already a set.\\

It follows from the above analysis that each of the maps $\Omega^{M,\lambda}_{i,B}$ making up $ \Map(\RuSpec(MA)(B), \Psi(B))$ 
corresponds to an inner commutative square in the diagram below for a given deformation $\delta: A \rarr MA$: 
\beq
\xymatrix{
	A \ar[ddr]_{\Lambda_{i,B}} \ar[dr]^{\delta} \ar[rrr] &&&id \ar@{=}[dl] \ar[ddl]^{\Gamma_B} \\
	& MA \ar[d]^{\Lambda_{i,B}^M} \ar[r] & id \ar[d]_{\Gamma_B}\\
	& B \ar[r] &\Psi(B)
} \nonumber
\eeq
thereby presenting the natural transformations $\Lambda_{i,B}^M$ as a factorization of $\Lambda_{i,B}$. By the functoriality of $\Psi$, one has:
\beq
\xymatrix{
	\Psi A \ar[dr] _{\Psi \delta} \ar[rrr] &&& \Psi id = id \ar@{=}[dl]  \ar[ddl]^{\Psi(\Gamma_B)}\\
	& \Psi(MA) \ar@{.>}[dr]_{\Psi( \Omega^{M,\lambda}_{i,B})} \ar[r] & \Psi id = id \ar[d] \\
	&& \Psi(\Psi(B))
} \nonumber
\eeq
where we used $\Omega^{M,\lambda}_{i,B}: MA \rarr id \rarr \Psi(B)$ from the previous paragraph, and $\Psi (id) = id$ since $\Psi$ is a morphism on $(\dX,\circ, id)$.\\

This is best presented in the form:
\beq
\xymatrixcolsep{3pc}
\xymatrix{
	\Psi A \ar[d]_{\Psi \delta} \ar[r] &id \ar[d]^{\Psi(\Gamma_B)} \\
	\Psi(MA) \ar[r]_{\Psi(\Omega^{M,\lambda}_{i,b})} & \Psi^2(B) 
} \nonumber
\eeq
which is a map $(\Psi A \Rarr \Psi(MA)) \rarr (id \Rarr \Psi^2(B))$, the first part of which would be a contribution to our naive definition of a cotangent set, and one would hope the converse holds, that is from the above diagram one retrieves the maps $\Omega_{i,B}^{M,\lambda}$, hence $\Der_{\Psi}(X,M)(B)$. Precisely because that is not necessarily the case however, we do not yet have an isomorphism of the form $\Der_F(X,M) \cong \bR \uh_s^{\LFx}$ as in \cite{TV4}. To circumvent this difficulty, we define the \textbf{cotangent set} of $\Psi$ at $x: X \equiv \RuSpec A \rarr \Psi$ in the direction of $M$ to be the object $\LPsix M$ such that:
\beq
\Der_{\Psi}(X,M)(B) \cong \Map(\LPsix M , \{id \Rarr \Psi^2(B)\}) \label{DerMapLPsiM}
\eeq
an isomorphism of sets, or $\Der_{\Psi}(X,M) \cong \Map(\LPsix M, \Psi^2)$, an isomorphism of functors $\dX \rarr \Set$, defined pointwise as in \eqref{DerMapLPsiM}. This being defined, the \textbf{cotanget set functor} of $\Psi$ at $x$ is given by:
\begin{align}
	\LPsix: \dX &\rarr \Set \nonumber \\
	M & \mapsto \LPsix M \nonumber
\end{align}

We define the \textbf{tangent set functor} of $\Psi$ at $x$ by $\TPsix = \Map(\LPsix, A )$, defined pointwise by $\TPsix M = \Map(\LPsix M, \{A = A\}) \in  \Set$. The latter expression gives the \textbf{tangent set} of $\Psi$ at $x$ in the direction $M$. This being defined, we finally have $\Cusp[\Psi, \eps] = \beeT_{\Psi, \eps}$, where $ \eps = A$, and we abused notation by writing $\beeT_{\Psi, \eps}$ for $\TPsix$, where $x: \RuSpec  (\eps)  \rarr \Psi$.\\

\subsection{Alternate simplicial formalism}
We defined $\Der_{\Psi}(X,M) = \Map(\RuSpec MA, \Psi)$, a $\Set$-valued functor, because we considered collections of maps $\Omega^{M, L}_B: \Lambda^M_B \rarr \Psi(B)$ for all $B \in \dX$. We now consider the cubical coverage of this formalism. Starting with the natural transformations $\Lambda^M_{i,B}: MA \Rarr B$, we consider the maximal cubic quiver $\Lambda^{\sigma,M}_B = (\Lambda^M_{i,B})$. Define:
\beq
\RuSpec^{\sigma} MA(B) = \maxMapgen(MA,B) \nonumber
\eeq
where $\maxMapgen(MA,B)$ indicates that we only consider the maximal quiver from $MA$ to $B$, and $gen$ is in reference to the fact that we work with generalized categories. With this definition, it then follows that we have:
\beq
\RuSpec^{\sigma} MA(B) = \Lambda^{\sigma,M}_B \nonumber 
\eeq
If now instead of considering select few morphisms $\Lambda^M_{i,B} \rarr \Gamma_B$ to populate a quiver we take a maximal quiver $\Omega^{\sigma,M}_{i,B} = ( \Lambda^M_{i,B} \rarr \Gamma_B)$ for each $i$, and if further we consider the full cube of all such maximal quivers, we get an object $\Omega^{\sigma, M}_B \in (\dX)^{gen}_{(MA,id,\Psi(B))} \in \Ho(\SetD)^{gen}$. We can then define:
\beq
\Der_{\Psi}^{\sigma}(X,M)(B) = \maxMapgen(\RuSpec^{\sigma} MA(B), \Gamma_B) = \Omega^{\sigma,M}_B \nonumber
\eeq
to be the maximal generalized morphism from $\RuSpec^{\sigma} MA(B)$ to $\Gamma_B$, an object of $ (\dX)^{gen}_{(MA,id,\Psi(B))} \in \Ho(\SetD)^{gen}$. Alternatively, note that for $\cC$ and $\cD$ categories, $\Fun(\cC,\cD^{gen}) \subset \Fun(\cC,\cD)^{gen}$, that is
\beq
\RuHom(\dX,(\dX)^{gen}) \subset \RuHom(\dX,\dX)^{gen} = \dtwoX^{gen} \label{dtwoXgen}
\eeq
This implies $\RuSpec^{\sigma} (MA) \in \dtwoX^{gen}$ and $\Psi \in \dtwoX^{gen}$, from which it follows $\Der_{\Psi}^{\sigma}(X,M) = \maxMapgen(\RuSpec^{\sigma} (MA), \Psi) \in (\dtwoX)^{gen}_{(\RuSpec^{\sigma}(MA),\Psi)} \in \Ho(\SetD)^{gen}$.\\

We define the \textbf{cotangent simplex} of $\Psi$ at $x: X \equiv \RuSpec^{\sigma} A \rarr \Psi$ in the direction of $M$ to be the object $\LPsix^{\sigma}M$ such that:
\beq
\Der_{\Psi}^{\sigma}(X,M)(B) \cong \maxMapgen(\LPsix^{\sigma} M , (id \Rarr \Psi^2(B))) \label{DerMapLPsiM2}
\eeq
an isomorphism in $\Ho(\SetD)^{gen}$, or $\Der_{\Psi}^{\sigma}(X,M) \cong \maxMapgen(\LPsix^{\sigma} M, \Psi^2)$, an isomorphism of functors $\dX \rarr \Ho(\SetD)^{gen}$ defined pointwise as in \eqref{DerMapLPsiM2}. To be precise, by the very definition of $\LPsix^{\sigma}$, we have $\LPsix^{\sigma} M \in (\dX)^{gen}_1$ so that $\maxMapgen(\LPsix^{\sigma}M, (id \Rarr \Psi^2(B))) \in (\dX)^{gen}_2$, hence $\maxMapgen(\LPsix^{\sigma} M,\Psi^2) \in \dtwoX^{gen}$ by \eqref{dtwoXgen}, so that the isomorphism $\Der_{\Psi}^{\sigma}(X,M) \cong \maxMapgen(\LPsix^{\sigma} M, \Psi^2)$ takes place in $\dtwoX^{gen}$.\\

Finally, we define the \textbf{tangent simplex} of $\Psi$ at $x$ by
\beq
\TPsix^{\sigma} = \maxMapgen(\LPsix^{\sigma}, A) \nonumber
\eeq
defined pointwise by $\TPsix^{\sigma} M = \maxMapgen(\LPsix^{\sigma} M, (A = A)) \in  \Ho(\SetD)^{gen}$. The latter expression gives the tangent simplex of $\Psi$ at $x$ in the direction $M$. This being defined, we finally have $\Cusp^{\sigma}[\Psi, \eps] = \beeT_{\Psi, \eps}^{\sigma}$, where $\eps = A$, and we abused notation by writing $\beeT_{\Psi, \eps}^{\sigma}$ for $\TPsix^{\sigma}$, with $x: \RuSpec^{\sigma} ( \eps) \rarr \Psi$.\\

\subsection{Connections between cubic and Set-based formalisms}
We can recover $\Der_{\Psi}(X,M)(B)$ for all $B \in \dX$ from its cubical counterpart by a truncation process: recall that if $L_i$ is an indexing set of flavors of maps $\Omega^{M,\lambda}_{i,B}: \Lambda^M_{i,B} \rarr \Gamma_B$, then one can define $\Omega^{\sigma,M,L_i}_{i,B} = ( \Omega^{M,\lambda}_{i,B})_{\lambda \in L_i}$. We have an obvious projection:
\beq
	\pi_{L_i}: \Omega^{\sigma,M}_{i,B} = (\Omega^{M,\lambda}_{i,B})_{\forall \lambda} \mapsto (\Omega^{M,\lambda}_{i,B})_{\lambda \in L_i} = \Omega^{\sigma, M,L_i}_{i,B} \nonumber 
\eeq
which we can generalize to:
\beq
\pi_L \Omega^{\sigma, M}_B = (\pi_{L_i}\Omega^{\sigma,M}_{i,B}) = (\Omega^{\sigma,M,L_i}_{i,B}) = \Omega^{\sigma,M,L}_B \nonumber
\eeq
as well as an iterative \textbf{deconstruction map} simply defined by:
\begin{align}
	\dagger (\Omega^{\sigma,M,L_i}_{i,B}) &= \{ \dagger \Omega^{\sigma,M,L_i}_{i,B} \}  \nonumber \\
	&= \{ \{ \Omega^{M,\lambda}_{i,B} \, | \, \lambda \in L_i \} \} \nonumber \\
	&= \{ \Omega^{M,L_i}_{i,B} \} \nonumber \\
	&= \Omega^{M,L}_B \nonumber
\end{align}
from which we have:
\begin{align}
	\dagger	\mathbf{\pi} \Der_{\Psi}^{\sigma}(X,M)(B) &= \dagger \pi \Omega^{\sigma,M}_B  \nonumber \\
	&:= \dagger \{ \pi_L \Omega^{\sigma,M}_L \}_L  \nonumber \\
	&=  \{ \dagger(\Omega^{\sigma,M,L_i}_{i,B}) \}\nonumber \\
	&=\{ \Omega^{M,L}_B \} \nonumber \\
	&=\Der_{\Psi}(X,M)(B) \nonumber 
\end{align}
Vice-versa, by definition $\Der_{\Psi}^{\sigma}(X,M)$ is constructed from $\Der_{\Psi}(X,M)$.\\

Using this truncation process, we now show those two notions of cusps we have defined, one being Set-based, the other one being simplicial, carry the same information. The same composite $\dagger \pi$ we used above also allows one to go from $\maxMapgen(\LPsix^{\sigma}M, (id \Rarr \Psi^2(B)))$ to $\Map(\LPsix M, (id \Rarr \Psi^2(B)))$:
\beq
\xymatrix{
	\Der_{\Psi}^{\sigma}(X,M)(B) \ar[d]_{\dagger \pi} \ar[r]^-{\cong} & \maxMapgen(\LPsix^{\sigma} M, (id \Rarr \Psi^2(B))) \ar@{.>>}[d]^{\dagger \pi} \\
	\Der_{\Psi}(X,M)(B) \ar[r]_-{\cong} & \Map(\LPsix M, (id \Rarr \Psi^2(B)))
} \nonumber
\eeq
from which it follows that:
\beq
\dagger \pi \LPsix^{\sigma}M = \LPsix M \nonumber
\eeq
Now:
\begin{align}
	\dagger \pi \TPsix^{\sigma} M &= \dagger \pi \maxMapgen(\LPsix^{\sigma} M, (A=A)) \nonumber \\
	&= \Map( \dagger \pi \LPsix^{\sigma} M, (A=A)) \nonumber \\
	&=\Map(\LPsix M, (A=A)) \nonumber \\
	&= \TPsix M \nonumber
\end{align}

\section{Strings}
Recall that objects of $\dAffCtildetau$ are stacks, which are modeled by functors $F: \CommC \rarr \SetD$. For $A \in \CommC$, $F(A)$ can be viewed as the manifestation, or realization, of $A$ via $F$. All those manifestations $F(A)$ of a phenomenon $F$ come together coherently by functoriality of $F$. In other terms, such a functor $F$ should therefore correspond to a physical system of some sort, since such systems vary in accordance to some laws, and this in a coherent fashion. Vice-versa, typical solutions to physical equations reflect certain laws, so assuming they depend functorially on them, if $D$ is a differential operator, $\Psi$ a solution to the equation $D \Psi = 0$, then one can identify $\Psi \equiv F: \CommC \rarr \SetD$.\\

Stacks are really objects of $\RHom(\dAffC^{\op}, \Top)$ with $\dAffC = L(\CommC^{\op})$, a Segal category, where $L$ denotes the Dwyer-Kan simplicial localization. Recall how this is performed (\cite{DK}): if $\cC$ is an ordinary category, $L\cC = F_*\cC[F_*W^{-1}]$ where $F_k\cC = F^{k+1}\cC$, $F\cC$ the free category on $\cC$ with generators $Fc$ for each non-identity morphism $c$ of $\cC$. Thus after a simplicial localization, equivalences have been sequentially contracted.\\

This will prove essential in our applications to Physical Mathematics in which we regard strings as modeling equivalences in $\CommC$, which are therefore no longer apparent in $\dAffC$, and since we will regard $\dAffCtildetau$, as done already in \cite{RG}, \cite{RG2}, \cite{RG4}, to model physical phenomena, this provides a reason for not being able to observe strings in nature.\\

Some physicists will argue that strings are fundamental objects. Already in Homotopy Type theory we see that this is not far-fetched. Indeed, one can regard sets as discrete homotopy types, thereby presenting homotopy types as fundamental objects. Those provide identification of elements through paths, which we regard as strings, in contrast to working with point-based Physics.\\

Practically for us then, working with strings corresponds to working with model categories, and after localization we are left with point-based Physics, which corresponds to the observable universe, modeled by stacks $F: \dAffC^{\op} \rarr \Top$.\\

To be precise, if $\cC$ is a simplicial, symmetric monoidal model category, assuming we have a model category structure on $\CommC$ as well in accordance with the assumptions put on a homotopical algebraic context, equivalences in $\CommC$ are defined on underlying simplicial objects, so $A \simeq A'$ in $\CommC$ means $\underline{A} \simeq \underline{A'}$ in $\SetD$. Those objects are therefore simplicially isomorphic in $\dAffC = L (\CommC^{\op})$, which means $\underline{A}_0 \cong \underline{A'}_0$, so we have an isomorphism between the objects of $A$ and $A'$. Define $a \in A_0$ and $a' \in A'_0$ in $\CommC$ to be equivalent: $a \sim a'$, if they are identified under the above isomorphism. We regard such an equivalence as a \textbf{string} between $a$ and $a'$ and we will sometimes write $a \stackrel{\gamma}{\sim} a'$. One important note at this point: those paths are formal, and there is no parametrization $\sigma$ along such a string, the way there would be in string theory. Also, there is no concept of propagation of such a string, hence no time parametrization $\tau$. In other terms, what we have is an algebraicization of physical string theory.\\

There are open and closed strings. Algebraically an open string is given by $a \sim b$ with $a \neq b$.\\
\beq
\setlength{\unitlength}{1cm}
\begin{picture}(2,2)(1,0)
	\put(0.5,1){$a$}
	\qbezier(1,1)(1.5,1.2)(2,1)
	\qbezier(2,1)(2.5,0.8)(3,1)
	\put(3.2,1){$b$}
	\put(1,1){\circle*{0.1}}
	\put(3,1){\circle*{0.1}}
\end{picture} \nonumber
\eeq
A closed string is given by $a \sim a$.\\
\beq
\setlength{\unitlength}{1cm}
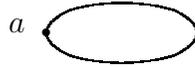
\begin{picture}(2,2)(1,0)
	\qbezier(1.3,1.3)(0.7,1)(1.3,0.7)
	\qbezier(1.3,0.7)(2,0.5)(2.7,0.7)
	\qbezier(2.7,0.7)(3.3,1)(2.7,1.3)
	\qbezier(2.7,1.3)(2,1.5)(1.3,1.3)
	\put(1,1){\circle*{0,1}}
	\put(0.5,1){$a$}
\end{picture} \nonumber
\eeq
In particular, every object of $A \in \CommC$ gives rise to its own closed string in a trivial way. Strings can split and join. For open strings, a join reads $a \sim b \wedge b \sim c \Rarr a \sim c$\\
\beq
\setlength{\unitlength}{1cm}
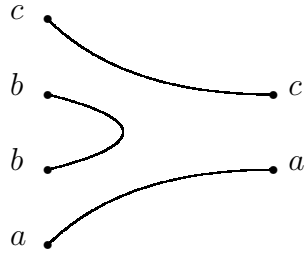
\begin{picture}(3,4)(0,-0.5)
	\put(-0.5,0){$a$}
	\put(-0.5,1){$b$}
	\put(-0.5,2){$b$}
	\put(-0.5,3){$c$}
	\put(0,0){\circle*{0.1}}
	\put(0,1){\circle*{0.1}}
	\put(0,2){\circle*{0.1}}
	\put(0,3){\circle*{0.1}}
	\qbezier(0,1)(2,1.5)(0,2)
	\qbezier(0,0)(1,1)(3,1)
	\qbezier(0,3)(1,2)(3,2)
	\put(3,1){\circle*{0.1}}
	\put(3.2,1){$a$}
	\put(3,2){\circle*{0.1}}
	\put(3.2,2){$c$}
\end{picture} \nonumber
\eeq
and a splitting reads $a \sim c \Rarr a \sim b \wedge b \sim c$, whenever $a \sim b$ or $b \sim c$:
\beq
\setlength{\unitlength}{1cm}
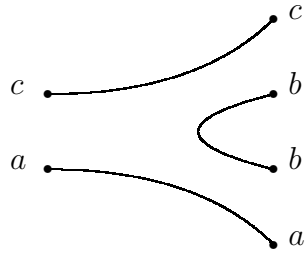
\begin{picture}(3,4)(0,-0.5)
	\put(3.2,0){$a$}
	\put(3.2,1){$b$}
	\put(3.2,2){$b$}
	\put(3.2,3){$c$}
	\put(3,0){\circle*{0.1}}
	\put(3,1){\circle*{0.1}}
	\put(3,2){\circle*{0.1}}
	\put(3,3){\circle*{0.1}}
	\qbezier(3,1)(1,1.5)(3,2)
	\qbezier(0,1)(2,1)(3,0)
	\qbezier(3,3)(2,2)(0,2)
	\put(0,1){\circle*{0.1}}
	\put(-0.5,1){$a$}
	\put(0,2){\circle*{0.1}}
	\put(-0.5,2){$c$}
\end{picture} \nonumber
\eeq

For closed strings this reads $(a \sim a) \wedge (a \sim a) \Rarr a \sim a$, a fold map
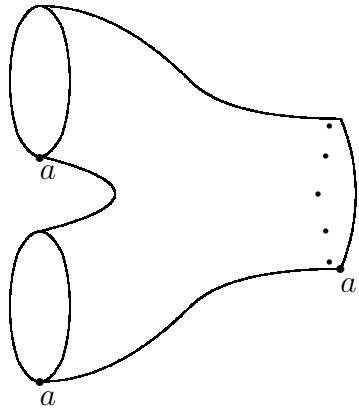
\beq
\setlength{\unitlength}{1cm}
\begin{picture}(5,6)(0,-0.5)
	\qbezier(0.7,4.7)(1,5.3)(1.3,4.7)
	\qbezier(1.3,4.7)(1.5,4)(1.3,3.3)
	\qbezier(0.7,3.3)(1,2.7)(1.3,3.3)
	\qbezier(0.7,4.7)(0.5,4)(0.7,3.3)
	\put(1,2.7){$a$}
	\put(1,2.98){\circle*{0.1}}
	\qbezier(0.7,1.7)(1,2.3)(1.3,1.7)
	\qbezier(1.3,1.7)(1.5,1)(1.3,0.3)
	\qbezier(0.7,0.3)(1,-0.3)(1.3,0.3)
	\qbezier(0.7,1.7)(0.5,1)(0.7,0.3)
	\put(1,-0.3){$a$}
	\put(1,0){\circle*{0.1}}
	\qbezier(5,3.5)(5.4,2.5)(5,1.5)
	\put(5,1.2){$a$}
	\put(5,1.5){\circle*{0.1}}
	\put(4.85,3.4){\circle*{0.06}}
	\put(4.8,3){\circle*{0.06}}
	\put(4.7,2.5){\circle*{0.06}}
	\put(4.8,2){\circle*{0.06}}
	\put(4.85,1.6){\circle*{0.06}}
	\qbezier(1,5)(2,5)(3,4)
	\qbezier(3,4)(3.5,3.5)(5,3.5)
	\qbezier(1,0)(2,0)(3,1)
	\qbezier(3,1)(3.5,1.5)(5,1.5)
	\qbezier(1,2)(3,2.5)(1,3)
\end{picture} \nonumber
\eeq
and $a \sim a \Rarr (a \sim a) \wedge (a \sim a)$, a diagonal map of sorts, with the reverse diagram. What is referred to as the world sheet in Physics for us corresponds to working in $\CommC$. A string in physical space-time is given by $X^{\mu}(\sigma, \tau)$, where $\mu$ is an index that keeps track of the dimension of space-time in which strings are propagating (\cite{GSW}). We argued a stack, modeled by a functor $F: \CommC \rarr \SetD$, corresponds to a physical field that satisfies some equations of motion, hence $F$ and $X^{\mu}$ are indistinguishable. If $A$ and $A'$ are two objects of $\CommC$, $A \simeq A'$, $a \in A$ and $a' \in A'$, identified after localization, we have a string between them: $a \sim a'$, and therefore a map between $Fa \in (FA)_0$ and $Fa' \in (FA')_0$. But $a \cong a'$ implies $Fa \cong Fa'$, hence $Fa \sim Fa'$ if we are willing to extend to $\SetD$ that notion of equivalence at the level of objects we had in $\CommC$. Then $Fa \stackrel{F\gamma}{\sim} Fa'$ plays the role of $X^{\mu}(\sigma, \tau)$, a string in $\SetD$. \\

Now one has to be mindful of the fact that $F\gamma$ is but a string, where $F \in \dAffCtildetau$, for $\gamma$ in $A \in \CommC$. $F$ corresponds to a given physical phenomenon, hence contains in its image all possible strings $F \gamma$ for $\gamma \in A$, for all laws $A \in \CommC$. It would be interesting to investigate how different string theories are encapsulated by different laws, and whether a functorial treatment such as the present one would shed light on M-theory, since a stack $F$ corresponds to a coherent unification of different string theories.\\

In contrast to string theory which has its own equations of motion, $F$ for us is just static, it tells us how a given phenonemon depends on physical laws $A \in \CommC$. Thus $\cX=\dAffCtildetau$ offers but a snapshot of a system at any moment. Dynamics as argued earlier corresponds to flows in $\RHomXX$. For simplicity, consider a flow built from a state $\eps$ of our system, an object of $\dX$, and some $\Psi \in \dtwoX$. $\Psi \eps$ is just another state of our system, and it is only $\Psi^{\infty} \eps $ that offers a flow of our system. In other terms the dynamics of $F$ in the present setting, and consequently that of the strings it is depicting, is given by $\Psi^{\infty}  \eps (F)$, which is formally defined by the following diagram, where we let $ \eps =  \eps_0 $ and $\Psi \eps_i =  \eps_{i+1} $:
\beq
\xymatrix{
	F \ar[r]^{\eps_0}  &\eps_0 F \ar[drrr] \ar@{|->}[r]^{\Psi} &  \eps_1  F \ar[drr] \ar@{|->}[r]^{\Psi} & \cdots \ar@{.>}[dr] \\
	&&&& \Psi^{\infty}  \eps  F 
} \nonumber
\eeq
\newline

To come back to strings, suppose $A \simeq A'$ in $\CommC$, so that $a \in A$, $a' \in A'$, being identified under the simplicial isomorphism $A \cong A'$ in $L \CommC$, are equivalent: $a \stackrel{\gamma}{\sim} a'$. If $F \in \dAffCtildetau$, we therefore have $Fa \stackrel{F\gamma}{\sim} Fa'$ in $\SetD$. Now let's consider a state $ \eps = \zeta$ of $\RHomXX$. In the same manner that a model for $F$ is given by a functor $F: \CommC \rarr \SetD$, a model for $\zeta$ is a functor $\zeta: \Mor(\CommC^{\op}, \SetD) \rarr \Mor(\CommC^{\op}, \SetD)$, that maps $F$ to $\zeta F$, and hence $F \gamma $ to $\zeta F \gamma$. Thus a state of our system does move strings around. Now for $\Psi \in \dtwoX$, $\Psi \zeta$ is some other state $ \phi$. Concretely, we have the following picture:
\beq
\xymatrix{
	F \gamma: Fa \sim Fa' \ar[d]_L \ar@{.>}[r]^{\zeta} &\zeta Fa \sim \zeta Fa' \ar@{.>}[r]^{\phi} & \phi \zeta Fa \sim \phi \zeta Fa' \\
	Fa \cong Fa' \ar[r]_{\zeta} & \zeta Fa \cong \zeta Fa' \ar[u]_{L^{-1}} \ar[r]_{\phi} & \phi \zeta Fa \cong \phi \zeta Fa' \ar[u]_{L^{-1}}
} \nonumber
\eeq
Regarding $\Psi$, the fashion in which it picks the state $\phi$ from a given state $\zeta$ of the system really depends on how it interprets natural laws, which are therefore inherently built into $\Psi$. Another interesting question is what mathematical formalism would allow us to make sense of such a concept.

\bigskip
\footnotesize
\noindent
\textit{e-mail address}: \texttt{rg.mathematics@gmail.com}.

\end{document}